\documentclass[12pt,twoside]{amsart}
\usepackage{amsthm,amsfonts,amssymb,amscd,euscript}

\usepackage[matrix,arrow]{xy}

\author{Florin Ambro} 
\address{RIMS, Kyoto University\\
Kyoto 606-8502, Japan.}
\email{ambro@kurims.kyoto-u.ac.jp}

\newcommand{\Q}{{\mathbb Q}}
\newcommand{\Z}{{\mathbb Z}}
\newcommand{\N}{{\mathbb N}}
\newcommand{\R}{{\mathbb R}}

\newcommand{\cO}{{\mathcal O}}
\newcommand{\cR}{{\mathcal R}}
\newcommand{\cS}{{\mathcal S}}

\newcommand{\bA}{{\mathbf A}}

\newcommand{\bD}{{\mathbf D}}
\newcommand{\bTheta}{{\mathbf \Theta}}
\newcommand{\bE}{{\mathbf E}}

\newcommand{\bK}{{\mathbf K}}

\newcommand{\bM}{{\mathbf M}}

\newcommand{\bR}{{\mathbf R}}

\newcommand{\Exc}{\operatorname{Exc}}
\newcommand{\Fix}{\operatorname{Fix}}
\newcommand{\bFix}{\operatorname{{\mathbf Fix}}}

\newcommand{\im}{\operatorname{Im}}

\newcommand{\mult}{\operatorname{mult}}

\newcommand{\Supp}{\operatorname{Supp}}

\theoremstyle{plain}
\newtheorem{thm}{Theorem}[section]

\newtheorem{lem}[thm]{Lemma}

\newtheorem{prop}[thm]{Proposition}

\theoremstyle{definition}

\newtheorem{rem}[thm]{Remark}

\newtheorem{ack}{Acknowledgments}

\theoremstyle{remark}

\begin{document}

\bibliographystyle{amsalpha+}
\title[Log canonical algebras]
{Restrictions of log canonical algebras of general type}
%\date{Aug 20, 2005}
\maketitle

\begin{abstract} We introduce a diophantine property of a
log canonical algebra, and use it to describe the restriction 
of a log canonical algebra of general type to a log canonical 
center of codimension one.
\end{abstract}

%\tableofcontents

%%%%%%%%%%%%%%%%%%%%%%%%%%%%%%%%%%%%%%%%%%%%%%%%%%%%%%%%%%%%%%%%%%%%%%
\setcounter{section}{-1}
%%%%%%%%%%%%%%%%%%%%%%%%%%%%%%%%%%%%%%%%%%%%%%%%%%%%%%%%%%%%%%%%%%%%%%
%%% Document name: lcalg.tex
%%% Last modified: Nov 22 2006
%%%%%%%%%%%%%%%%%%%%%%%%%%%%%%%%%%%%%%%%%%%%%%%%%%%%%%%%%%%%%%%%%%%%%%

%%%%%%%%%%%%%%%%%%%%%%%%%%%%%%%%%%%%%%%
%%%%%%%%%%%%%%%%%%%%%%%%%%%%%%%%%%%%%%%

\section{Introduction}

%%%%%%%%%%%%%%%%%%%%%%%%%%%%%%%%%%%%%%%
%%%%%%%%%%%%%%%%%%%%%%%%%%%%%%%%%%%%%%%
\footnotetext[1]{This work was supported by a 
Twenty-First Century COE Kyoto Mathematics Fellowship,
and by the JSPS Grant-in-Aid No 17740011.
}
\footnotetext[2]{2000 Mathematics Subject Classification. 
Primary: 14C20. Secondary: 14E30.}

In an inductive study of log canonical algebras, it 
is important to understand their restrictions to log 
canonical centers. Let $X$ be a nonsingular complex variety,
$Y\subset X$ a nonsingular divisor and $\pi\colon X\to S$ a 
projective morphism.
The adjunction formula $(K+Y)\vert_Y=K_Y$ induces a 
homomorphism of graded $\cO_S$-algebras
$$
\cR_{X/S}(K+Y)\to \cR_{Y/S}(K_Y),
$$
where 
$
\cR_{X/S}(K+Y)=\bigoplus_{i=0}^\infty \pi_*\cO_X(iK+iY)
$ 
is the log canonical algebra of $(X/S,Y)$ and 
$
\cR_{Y/S}(K_Y)=\bigoplus_{i=0}^\infty \pi_*\cO_Y(iK_Y)
$ 
is the log canonical algebra of $Y/S$. The image of this 
homomorphism is a graded subalgebra, denoted
\begin{equation}\label{inc}
\cR_{X/S}(K+Y)\vert_Y\subseteq \cR_{Y/S}(K_Y).
\end{equation}
Siu's~\cite{Siu98} invariance of plurigenera of varieties
of general type can be restated as follows: if $S$ is a 
smooth curve and $Y$ is a smooth $\pi$-fiber of 
general type, the restricted algebra $\cR_{X/S}(K+Y)\vert_Y$
coincides with $\cR_{Y/S}(K_Y)$.
Kawamata~\cite{Kaw99can, Ka99} and Nakayama~\cite{Nak98, Nak04}
obtained singular versions of this result, and extended it 
in a different direction: if $K+Y$ is $\pi$-big and its 
relative Iitaka map maps $Y$ birationally onto its image, 
the restricted algebra $\cR_{X/S}(K+Y)\vert_Y$ coincides 
with $\cR_{Y/S}(K_Y)$ in degrees $i\ge 2$. 

A characterization of the restricted algebra was also expected 
in the logarithmic case (see Nakayama~\cite[Theorem 4.9]{Nak98}), 
but a new point of view was necessary, since it was known that 
the inclusion in~\eqref{inc} may be strict in all degrees in
this case.
The new idea, due to Hacon and M\textsuperscript{c}Kernan, 
is that the restricted algebra is equivalent with the log 
canonical algebra of a log structure defined not necessarily
on $Y$, but on a birational model of $Y$. 
Hacon and M\textsuperscript{c}Kernan (\cite{JH05}, Theorem 4.3) 
obtained this description in the birational log Fano set-up of 
a prelimiting flip and assuming the validity of the log Minimal 
Model Program in smaller dimensions. In this paper we 
sharpen their result, and extend it to log varieties of general 
type. The other contribution of this paper is a new diophantine 
property of a log canonical algebra (Lemma~\ref{k+}), which
is based on ideas of Shokurov and Viehweg. We hope these two
new tools will be useful for bounding log canonical models 
(see Koll\'ar~\cite{Kol92} for basic open problems). 
We also present two applications in \S 4.

Our main result is as follows (we work over an algebraically 
closed field of characteristic zero).

\begin{thm}\label{mr} Let $X$ be a nonsingular algebraic variety
and $B$ an $\R$-divisor on $X$ such that $\Supp(B)$ is 
a simple normal crossings divisor. Assume that $Y$ is a
component of $B$ with $\mult_Y(B)=1,\lfloor B-Y\rfloor=0$.
Denote $B_Y=(B-Y)\vert_Y$, so that by adjunction we have
$
(K+B)\vert_Y=K_Y+B_Y.
$
Let $\pi\colon X\to S$ be a projective surjective morphism,
and assume
\begin{itemize}
\item[(a)] $K+B\sim_\Q A+C$, where $A$ is a $\pi$-ample 
$\R$-divisor and $C$ is an effective $\R$-divisor with 
$\mult_Y(A)=\mult_Y(C)=0$.
\item[(b)] $(Y,B_Y)$ has canonical singularities 
in codimension at least two.
\end{itemize}
Define 
$
\Theta=\max(B_Y-\lim_{i\to \infty}
\frac{(\bFix(iK+iB)\vert_Y)_Y}{i},0)
$
where the maximum is taken componentwise. For every $n\ge 1$ we 
have natural inclusions
$$
\im(\pi_*\cO_X(nK+nB)\to \pi_*\cO_Y(nK_Y+nB_Y))
\subseteq \pi_*\cO_Y(nK_Y+n\Theta).
$$
The following properties hold:
\begin{itemize}
\item[(1)]
The inclusion is an equality if $n\ge 2$ and $\{nB\}\le B$,
or if $n=1$ and $\pi(Y)\ne \pi(X)$. 
\item[(2)] Assume that $B$ has rational coefficients, and
the log canonical divisor $K_Y+\Theta$ has a Zariski 
decomposition relative to $S$. Then $\Theta$ has rational 
coefficients, and the graded $\cO_S$-algebra
$\bigoplus_{n=0}^\infty \pi_*\cO_Y(nK_Y+n\Theta)$ is finitely 
generated.
\end{itemize}
\end{thm}

In the definition of $\Theta$, $(\bFix(iK+iB)\vert_Y)_Y$
is the trace on $Y$ of the restriction to $Y$ of the fixed 
$\R$-b-divisor of $\pi_*\cO_X(iK+iB)$. Precisely, 
assume that $i(K+B)$ is relatively mobile at $Y$
and choose a birational modification $\mu\colon X'\to X$
such that the mobile part $M_i$ of $\mu^*(iK+iB)$ is 
relatively free, and the proper transform $Y'$ of $Y$
on $X$ is normal. Then $(\bFix(iK+iB)\vert_Y)_Y$ is the
push forward of $(\mu^*(iK+iB)-M_i)\vert_{Y'}$ via 
the birational map $Y'\to Y$.

 The assumption (b) is necessary for the restricted 
algebra $\cR_{X/S}(K+B)\vert_Y$ to have a presentation 
as a log canonical algebra on $Y$. A similar result holds 
when $(Y,B_Y)$ has only Kawamata log terminal singularities 
(Theorem~\ref{mth}), but one has to pass to a birational
model of $Y$ so that $\Theta$ takes into account all 
valuations of $Y$ whose log discrepancy with respect to 
$(Y,B_Y)$ is less than one (see Section 2 for details). 

As for the proof of Theorem~\ref{mr}, we recommend that the 
reader first consults Lemma~\ref{cd}, for an argument
modulo the log Minimal Model Program in the same dimension. 
The proof of (1) is based on Siu's idea, with modifications 
by Kawamata and Nakayama. Siu's method~\cite{Siu98} of dealing
with pluricanonical sections is to 
view $nK_Y$ as $K_Y+(n-1)K_Y$, and to pass by induction a 
property from $(n-1)K_Y$ to $nK_Y$. This still works in 
the logarithmic case, provided we replace the given 
boundary by a canonical sequence of boundaries, satisfying 
certain arithmetic properties (Lemmas~\ref{Sep22} and
~\ref{k+}). 
The real coefficients of the boundary pose no problem, 
since Kawamata-Viehweg vanishing is known to hold in 
the real case. Also, by diophantine approximation, there 
are infinitely many positive integers $n$ satisfying 
the inequality $\{nB\}\le B$. 

The proof of (2) is based on a new diophantine property of 
a log canonical algebra (Lemma 1.5), 
and a criterion for a real nef divisor to be rational and 
semiample~\cite{Semcr}. The former is a combination of 
Shokurov's~\cite{Sh03} ideas on diophantine properties of graded 
algebras and Viehweg's~\cite{Viehweg82} method of dealing with 
pluricanonical sections in his proof of the weak positivity 
of the push forwards of relative pluricanonical sheaves, and the 
latter generalizes Kawamata's criterion~\cite{Kaw87} that 
log canonical rings of
general type are finitely generated if Zariski decomposition
exists. Conversely, the log Minimal Model Program (with real 
boundaries) in the dimension of $\tilde{Y}$ implies 
the existence of Zariski decomposition for the big log 
canonical divisor $K_{\tilde{Y}}+\Theta$ in Theorem~\ref{mr}
(see Lemma~\ref{zd}).

 The reader will notice that we make heavy use of Shokurov's 
new terminology of b-divisors, instead of multiplier ideal sheaves,
which are common in this context. The logarithmic implementation 
of Siu's idea involves taking a log canonical divisor 
out of the round-up in Lemma~\ref{k+}, on sufficiently 
high birational models of a given variety, and since the 
log canonical divisor does not have integer coefficients, we 
work directly with divisors on these high models. B-divisors 
are a very useful notation for making computations on these high
models, finitely many at a time, but the reader may avoid them by 
simply introducing 
notation for these models. We are unable to encode 
this argument on the base variety, in terms of multiplier ideal 
sheaves. The reader interested in this may consult the  
arguments of Takayama~\cite[Theorem 4.1]{Takayama06} and
Hacon-M\textsuperscript{c}Kernan~\cite[Corollary 3.17]{HJ05}.
 
 Finally, we expect that the hypothesis $(a)$ in Theorem~\ref{mr} 
can be weakened to $(a')$: $K+B\sim_\Q C$, where $C$ is effective, 
$\pi$-big and $\mult_Y(C)=0$. This may follow from an 
algebrization of the method introduced by Siu for the 
invariance of plurigenera of manifolds of non-general type 
(see~\cite{Siu05}).

\begin{ack} I would like to thank Professors Yujiro
Kawamata, James M\textsuperscript{c}Kernan and 
Noboru Nakayama for useful discussions.
\end{ack}

%%%%%%%%%%%%%%%%%%%%%%%%%%%%%%%%%%%%%%%
%%%%%%%%%%%%%%%%%%%%%%%%%%%%%%%%%%%%%%%

\section{Preliminary}

%%%%%%%%%%%%%%%%%%%%%%%%%%%%%%%%%%%%%%%
%%%%%%%%%%%%%%%%%%%%%%%%%%%%%%%%%%%%%%%

%%%%%%%%%%%%%%%%%%%%%%%%%%%%%%%%%%%%%%%
%%%%%%%%%%%%%%%%%%%%%%%%%%%%%%%%%%%%%%%

\subsection{Boundary arithmetic}

%%%%%%%%%%%%%%%%%%%%%%%%%%%%%%%%%%%%%%%
%%%%%%%%%%%%%%%%%%%%%%%%%%%%%%%%%%%%%%%
\begin{lem}\label{Sep22} 
For $b\in [0,1]$ and $n\ge 1$ define
$
b_n=\max(b,\frac{1}{n}\lceil (n-1)b\rceil).
$
The following properties hold:
\begin{itemize}
\item[(1)] $b=b_n$ if and only if $\{nb\}\le b$.
\item[(2)] $b\le b_n\le b+\frac{1}{n}$.
\item[(3)] $\lfloor b\rfloor+\lceil (n-1) 
b_{n-1}\rceil\le n b_n$ for $n\ge 1$.
\end{itemize}
\end{lem}

\begin{proof} Properties (1) and (2) are easy to see. 
Property (3) is clear if $b=0$ or $1$, so let $b\in (0,1)$. 
The claim is then equivalent to 
$$
(n-1)b \le \lfloor n b_n\rfloor.
$$
Assume first that $\{nb\}\le b$, that is $b_n=b$. Then 
$\lfloor n b\rfloor-(n-1)b=b-\{nb\}\ge 0$.
Assume now that $\{nb\}\ge b$, that is 
$nb_n=\lceil (n-1)b\rceil$.
Then $\lceil (n-1)b\rceil-(n-1)b=\{(n-1)b\}\ge 0$.
\end{proof}

\begin{lem}\label{inq} 
Let $n$ be a positive integer and $b,e\in \R_{\ge 0}$ 
such that $e-bn\in \Z$. Then $\lceil -b+\frac{e}{n}\rceil\le e$.
\end{lem}

\begin{proof} Let $e-bn=p\in \Z$. For $p\ge 0$, we have
$$
e-\lceil -b+\frac{e}{n}\rceil=p-\lceil \frac{p}{n}\rceil+bn\ge 0.
$$
If $p\le 0$, then 
$\lceil -b+\frac{e}{n}\rceil=\lceil \frac{p}{n}\rceil\le 0
\le e$.
\end{proof}

\begin{lem}\label{div} 
Let $a,c,d,\gamma\in \R$ such that $a>-1$, $a-\gamma c\ge -1$ and
$\gamma>0$. If $n$ is an integer such that $n \ge 1+\gamma^{-1}$, then
$$
\lfloor \frac{\lceil a+c+nd\rceil}{n}\rfloor\le\lceil a+d\rceil.
$$
\end{lem}

\begin{proof} Since $a+1+c\le (a+1)(1+\gamma^{-1})$ and 
$\lceil a+1+d\rceil-d\ge a+1$, we obtain
$1+\gamma^{-1}\ge \frac{a+1+c}{\lceil a+1+d\rceil-d}$.
Therefore
$$
n\ge \frac{a+1+c}{\lceil a+1+d\rceil-d}.
$$
This is equivalent to the conclusion, by a straightforward
computation.
\end{proof}

%%%%%%%%%%%%%%%%%%%%%%%%%%%%%%%%%%%%%%%
%%%%%%%%%%%%%%%%%%%%%%%%%%%%%%%%%%%%%%%

\subsection{B-divisors, log pairs, log varieties}

%%%%%%%%%%%%%%%%%%%%%%%%%%%%%%%%%%%%%%%
%%%%%%%%%%%%%%%%%%%%%%%%%%%%%%%%%%%%%%%

We refer the reader to~\cite{Semcr} for standard
definitions on Shokurov's b-divisors, log pairs and log
varieties. Just to fix the notation, recall that 
a log pair $(X,B)$ is a normal complex variety $X$
endowed with an $\R$-divisor $B$ such that $K+B$
is $\R$-Cartier. A log variety is a log pair whose
boundary $B$ is effective. 
The discrepancy $\R$-b-divisor of a log pair $(X,B)$ is
$$
\bA(X,B)=\bK-\overline{K+B},
$$
where $\bK$ is the canonical b-divisor of $X$ and
$\overline{K+B}$ is the Cartier closure of the log
canonical class. If $(X,B)$ has log canonical
singularities, let $\bR$ be the reduced b-divisor
of all prime b-divisors of $X$ which have zero log
discrepancy with respect to $(X,B)$. Define
$
\bA(X,B)^*=\bA(X,B)+\bR,
$
so that $\lceil \bA(X,B)^*\rceil\ge 0$.

Let $\pi\colon X\to S$ be a proper morphism from 
a normal variety $X$ and let $D$ be an $\R$-Cartier
$\R$-divisor on $X$. We denote by $\bM(D)$ the mobile
b-divisor of $D$ relative to $S$, $\bFix(D)=
\overline{D}-\bM(D)$, $\bD_i(D)=\frac{1}{i}\bM(iD)$
and 
$$
\cR_{X/S}(D)=\bigoplus_{i=0}^\infty \pi_*\cO_X(iD).
$$
Note that $\Fix(D)$, the trace of $\bFix(D)$ on $X$,
is the fixed part of $D$ relative to $S$ in the usual 
sense. Locally over $S$,
$$
\Fix(D)=\inf\{(a)+D; a\in \pi_*\cO_X(D)\setminus 0\}.
$$

\begin{lem}[Terminal resolution]\label{fn}
Let $(X,B)$ be a log pair with Kawamata log terminal 
singularities. Then the set of prime b-divisors $E$ 
of $X$, having log discrepancy $a(E;X,B)$ less than $1$, 
is finite.
\end{lem}

\begin{proof} We may assume that $X$ is smooth and
$\Supp(B)$ has simple normal crossings. Consider the
set of pairs of distinct prime divisors on $X$
$$
\cS=\{(E_1,E_2); E_1\cap E_2\ne \emptyset,
a(E_1;X,B)+a(E_2;X,B)\le 1\}.
$$
If $\cS$ is empty, it is easy to see that $a(E;X,B)> 1$ 
for every prime b-divisor $E$ which is exceptional over
$X$. If $\cS$ is nonempty, define
$$
a=\min_{(E_1,E_2)\in \cS}\min(a(E_1;X,B),a(E_2;X,B)).
$$
Since $(X,B)$ has Kawamata log terminal singularities,
we have $a>0$ and 
$$
a(E_1;X,B)+a(E_2;X,B)\ge 2a \mbox{ for }(E_1,E_2)\in \cS.
$$
Let $X_1\to X$ be the composition of the 
blow-ups of $X$ in $E_1\cap E_2$, after all 
$(E_1,E_2)\in \cS$, and let $\mu_1^*(K+B)=K_{X_1}+B_{X_1}$.
Then $X_1$ is smooth and $\Supp(B_{X_1})$ has simple normal 
crossings. Consider the set of distinct 
prime divisors on $X_1$
$$
\cS_1=\{(E_1,E_2); E_1\cap E_2\ne \emptyset,
a(E_1;X_1,B_{X_1})+a(E_2;X_1,B_{X_1})\le 1\}.
$$
If $\cS_1$ is empty, we are done. 
Otherwise, 
$$a(E_1;X_1,B_{X_1})+a(E_2;X_1,B_{X_1})\ge 3a \mbox{ for } (E_1,E_2)
\in \cS_1,
$$ 
and we repeat the process: let $X_2\to X_1$ be the composition 
of the blow-ups of $X_1$ in $E_1\cap E_2$, after all 
$(E_1,E_2)\in \cS$, and so on. After $n$ blow-ups, either 
$\cS_n$ is empty, or 
$$
a(E_1;X_n,B_{X_n})+a(E_2;X_n,B_{X_n})\ge (n+1)a \mbox{ for } 
(E_1,E_2)\in \cS_n.
$$
Therefore there exists $n\le \lceil a^{-1}\rceil$ such that
$\cS_n=\emptyset$, that is $(X_n,B_{X_n})$ has terminal 
singularities in codimension at least two.
\end{proof}

The next lemma is our key tool, a diophantine 
property of log canonical algebras. Inspired
by Shokurov's notion of asymptotic saturation of a graded 
algebra~\cite{Sh03}, it is the logarithmic version of
Viehweg's~\cite{Viehweg82} method of dealing with 
pluricanonical sections, which is in fact similar to 
Siu's~\cite{Siu98}. 
If $X$ is nonsingular and $B=0$,
property (2) is implicit in Viehweg's proof of the weak 
positivity of the push forwards of relative pluricanonical 
sheaves. If $B=0$, Lemma~\ref{k+} is a 
priori stronger than the asymptotic saturation property 
of $\cR_{X/S}(K_X)$ with respect to the log variety $X$, 
but the two notions coincide if $\cR_{X/S}(K_X)$ 
is finitely generated. In general, the two notions differ.
Consider for example 
$\cR_{{\mathbb C}/{\mathbb C}}(K_{\mathbb C}+b\cdot 0)$,
where $b\in (0,1)$ and ${\mathbb C}/{\mathbb C}$ is
the identity map. Then Lemma~\ref{k+} says that
$\lfloor nb\rfloor\le \lceil (n-1)b\rceil$ for every
$n\ge 1$, with equality if $\{nb\}\le b$. On the other
hand, $\cR_{{\mathbb C}/{\mathbb C}}(K_{\mathbb C}+b\cdot 0)$
is asymptotically saturated with respect to the log variety 
$({\mathbb C},b\cdot 0)$ if and only if 
$\lceil (n-1)b\rceil\le \lfloor nb\rfloor$ for sufficiently
large and divisible integers $n$, that is $b\in \Q$. 
The reader may consult~\cite{Semcr,ToricFga} for more
on asymptotically saturated graded algebras. Finally,
Lemma~\ref{k+} has an analogue for anti-log canonical 
algebras, and it seems to be peculiar to these two type 
of algebras.

\begin{lem}\label{k+} Let $(X,B)$ be a log variety with
log canonical singularities, and $\pi\colon X\to
S$ a proper morphism. Let $\bR$ be the reduced
b-divisor of all prime b-divisors of $X$ which have zero
log discrepancy with respect to $(X,B)$. Note that $\bR=0$
if and only if $(X,B)$ has Kawamata log terminal 
singularities. 
Let $n$ be a positive integer such that 
$\pi_*\cO_X(nK+nB)\ne 0$. Then the following properties
hold:
\begin{itemize}
\item[(1)] For every $i\in n\N$ we have an inclusion
$$
\pi_*\cO_X(nK+nB)\subseteq \pi_*\cO_X(\lceil\bK+\bR+
(n-1)\bD_i(K+B)\rceil).
$$
The sheaf on the right-hand side is independent
of $i$ sufficiently large and divisible.
\item[(2)] Equality holds in (1) if $i=n$, or if
$\{nB\}\le B$ and $i\in n\N$.
\end{itemize}
\end{lem}

\begin{proof} (1) This follows from the case $i=n$ of (2),
since $\bD_n(K+B)\le \bD_i(K+B)$ for $n\vert i$.

(2a) We show equality holds for $i=n$.
For the direct inclusion, let $a\in k(X)^\times$ with 
$(a)+nK+nB\ge 0$. In particular, 
$\overline{(a)}+n\bD_n(K+B)\ge 0$. 
Since $\bD_n(K+B)\le \overline{K+B}$, we have
$$
\bK+\bR+(n-1)+n\bD_n(K+B)\ge \bA(X,B)^*+n\bD_n(K+B).
$$
Since $\lceil \bA(X,B)^*\rceil\ge 0$, we obtain
$$
\overline{(a)}+\lceil\bK+\bR+(n-1)\bD_n(K+B)\rceil
\ge \lceil \bA(X,B)^*+(\overline{(a)}+n\bD_n(K+B))\rceil\ge 0.
$$
We now consider the opposite inclusion. By Hironaka, there 
exists a proper birational morphism $\mu\colon X'\to X$ and 
a Cartier divisor $M_n$ on $X'$ such that
\begin{itemize}
\item[(i)] $M_n\le \mu^*(nK+nB)$ is $(\pi\circ\mu)$-free.
\item[(ii)] The inclusion $(\pi\circ\mu)_*\cO_{X'}(M_n)
\subset \pi_*\cO_X(nK+nB)$ is an equality.
\item[(iii)] $X'$ is nonsingular, $\Supp(M_n)\cup\Supp(B_{X'})$ 
is a simple normal crossings divisor, where 
$\mu^*(K+B)=K_{X'}+B_{X'}$ is the log pullback.
\end{itemize}
Thus $\bD_n(K+B)$ is the Cartier closure of $\frac{1}{n}M_n$. 
Let $R=\sum_{\mult_E(B_{X'})=1}E$. By construction, 
$F_n=\mu^*(nK+nB)-M_n$ is an effective $\R$-divisor. We have
$$
\lceil K_{X'}+R+\frac{n-1}{n}M_n\rceil=
M_n+\lceil-B_{X'}+R+\frac{1}{n}F_n\rceil.
$$
Since $(X,B)$ has log canonical singularities and
$B$ is effective, the divisor $\lceil -B_Y+R\rceil$ is 
effective and $\mu$-exceptional. Since $F_n$ is effective, 
we obtain
$$
0\le \lceil-B_{X'}+R+\frac{1}{n}F_n\rceil.
$$
On the other hand, $F_n-nB_{X'}=nK_{X'}-M_n$ is an integral divisor.
By Lemma~\ref{inq}, we obtain
$$
\lceil -B_{X'}+R+\frac{1}{n}F_n\rceil\le F_n+G_n,
$$
where 
$G_n=\sum_{\mult_P(B_{X'})<0}\lceil \mult_P(-B_{X'}+
\frac{1}{n}F_n)\rceil \cdot P$, where the sum runs after 
the prime divisors in $X'$ where $B_{X'}$ is negative. 
Combining the above inequalities, we obtain
$$
M_n\le \lceil K_{X'}+R+\frac{n-1}{n}M_n\rceil \le
\mu^*(nK+nB)+G_n.
$$
The Cartier divisor $G_n$ is effective and 
$\mu$-exceptional, since $B$ is effective. In particular,
$$
\mu_*\cO_{X'}(\mu^*(nK+nB)+G_n)=\cO_X(nK+nB).
$$
Therefore we obtain inclusions
$$
\mu_*\cO_{X'}(M_n)\subseteq 
\mu_*\cO_{X'}(\lceil K_{X'}+R+\frac{n-1}{n}M_n\rceil)
\subseteq \cO_X(nK+nB).
$$
Since $\pi_*\mu_*\cO_{X'}(M_n)=\pi_*\cO_X(nK+nB)$, 
we obtain 
$$
\pi_*\cO_X(nK+nB)=
\pi_*\mu_*\cO_{X'}(\lceil K_{X'}+R+\frac{n-1}{n}M_n\rceil).
$$
(2b) Assume now that $\{nB\}\le B$ and $i\in n\N$.
We have inclusions
\begin{align*}
\pi_*\cO_X(\lceil \bK+\bR+(n-1)\bD_i(K+B)\rceil) & 
 \subseteq \pi_*\cO_X(\lceil \bK+\bR+(n-1)\overline{K+B}\rceil)\\
  & \subseteq \pi_*\cO_X(\lceil nK+\lfloor B\rfloor+
(n-1)B\rceil) \\
  & \subseteq \pi_*\cO_X(nK+nB),
\end{align*}
where the first inclusion holds by 
$\bD_i(K+B)\le \overline{K+B}$, and the third by
Lemma~\ref{Sep22}.
By (2a), all inclusions are equalities.
\end{proof}

%%%%%%%%%%%%%%%%%%%%%%%%%%%%%%%%%%%%%%%
%%%%%%%%%%%%%%%%%%%%%%%%%%%%%%%%%%%%%%%

\subsection{Zariski decomposition}

%%%%%%%%%%%%%%%%%%%%%%%%%%%%%%%%%%%%%%%
%%%%%%%%%%%%%%%%%%%%%%%%%%%%%%%%%%%%%%%
We refer the reader to Nakayama~\cite{Nak04} for an 
excellent introduction to the Zariski decomposition 
problem.  
Several higher dimensional versions of the two-dimensional 
Zariski decomposition have been proposed, but they all 
coincide for big divisors.

Consider a proper surjective morphism $\pi\colon X\to S$
and a $\pi$-big $\R$-divisor $D$ on $X$. We say 
that $D$ has a Zariski decomposition, 
relative to $S$, if there exists a birational contraction 
$\mu\colon X'\to X$ and a $(\pi\circ\mu)$-nef and 
$(\pi\circ\mu)$-big $\R$-divisor $P$ on $X'$ such that
\begin{itemize}
\item[(i)] $P\le \mu^*D$;
\item[(ii)] $\cR_{X'/S}(P)=\cR_{X/S}(D)$.
\end{itemize}

\begin{lem}\label{zd} Let $(X,B)$ be a log variety with log
canonical singularities, and $\pi\colon X\to
S$ be a proper surjective morphism such that
$K+B$ is $\pi$-big. Assume the log Minimal Model Program
(with real boundaries) is valid in dimension $\dim(X)$.

Then $K+B$ has a Zariski decomposition,
and $P$ is the pullback to a suitable model of the log
canonical class of the log minimal model of $(X,B)$.
\end{lem}

\begin{proof} If log Minimal Model Program holds for
$(X/S,B)$, we obtain a birational map to a log minimal
model over $S$
$$
\Phi\colon (X,B)\dashrightarrow (X',B').
$$
If we consider the normalization of the graph of $\Phi$, 
this means that we have a Hironaka hut
 \[ \xymatrix{
& X'' \ar[dl]_\mu \ar[dr]^{\mu'} &  \\
X \ar@{.>}[rr]^\Phi & & X'
 } \]
with the following properties:
\begin{itemize}
\item[(i)] $(X',B')$ is a log variety with log canonical
singularities and $K_{X'}+B'$ is relatively nef and big;
\item[(ii)] $F=\mu^*(K+B)-{\mu'}^*(K_{X'}+B')$ 
is effective and $\mu'$-exceptional.
\end{itemize}
Denote $P={\mu'}^*(K_{X'}+B')$. By (ii), we have
$\cR_{X/S}(K+B)=\cR_{X''/S}(P)$. Therefore 
$P\le \mu^*(K+B)$ is a Zariski decomposition.
\end{proof}

\begin{rem}[Real logMMP] The largest category 
in which logMMP is expected to work is that of relative log
varieties $(X/S,B)$ with log canonical singularities.
The Cone and Contraction Theorems are known (\cite{KMM,Qlv}). 
The existence of extremal flips is numerical, hence it 
follows from the existence of flips with rational boundary. 
The termination of log flips is known in dimension $3$ 
(Shokurov~\cite{Sh96}) and is open in dimension at least 
$4$ (see Shokurov~\cite{Sh04} for more on termination).
\end{rem}

%%%%%%%%%%%%%%%%%%%%%%%%%%%%%%%%%%%%%%%
%%%%%%%%%%%%%%%%%%%%%%%%%%%%%%%%%%%%%%%

\section{The boundary of the induced log canonical algebra}

%%%%%%%%%%%%%%%%%%%%%%%%%%%%%%%%%%%%%%%
%%%%%%%%%%%%%%%%%%%%%%%%%%%%%%%%%%%%%%%

Let $(X,B)$ be a log pair with log canonical singularities 
and let $\pi\colon X\to S$ be a proper morphism. Assume that 
$Y$ is a normal prime divisor in $X$ with $\mult_Y(B)=1$,
and there exists a positive integer $l$ such that $lK+lB$ is 
$\pi$-mobile at $Y$.

By codimension one adjunction, there exists a 
canonically defined log pair structure
$(Y,B_Y)$ on $Y$ such that 
$$
(K+B)\vert_Y=K_Y+B_Y.
$$
For every $i\in l\N$, let $\bFix(iK+iB)$ be the 
fixed $\R$-b-divisor of $iK+iB$ relative to $S$. By
assumption, $\bFix(iK+iB)$ is b-$\R$-Cartier
and it has multiplicity zero at $Y$. Therefore
its restriction $\bFix(iK+iB)\vert_Y$ is a well
defined b-$\R$-Cartier $\R$-b-divisor of $Y$. 
Define
$$
\bTheta=\sum_{E}
\max(1-a(E;Y,B_Y)-\lim_{i\to \infty}
\mult_E(\bFix(iK+iB)\vert_Y),0)E,
$$
where the sum runs after all prime b-divisors
$E$ of $Y$ and $a(E;Y,B_Y)$ denotes the log
discrepancy of $E$ with respect to $(Y,B_Y)$.
It is clear that $\bTheta$ is supported
by the prime b-divisors $E$ with
$a(E;Y,B_Y)\in [0,1)$. In particular, $\bTheta$
is a well defined effective $\R$-b-divisor of $Y$.

\begin{lem}\label{cut} Let $\mu\colon X'\to X$ be a 
birational contraction and let $\mu^*(K+B)=K_{X'}+B_{X'}$ 
be the log pullback. 
\begin{itemize}
\item[(i)] The log pairs $(X,B)$ and $(X',B_{X'})$ induce
the same $\bTheta$.
\item[(ii)] Assume that $B\ge 0$ and $K_{X'}+B'$ is 
$\R$-Cartier, where $B'=\max(B_{X'},0)$. Then the log
varieties $(X,B)$ and $(X',B')$ induce the same $\bTheta$.
\end{itemize}
\end{lem}

\begin{proof} Property (i) is clear. For (ii), we have
$$
K_{X'}+B'=\mu^*(K+B)+F,
$$
where $F$ is effective and $\mu$-exceptional. In
particular,
$$
\pi_*\cO_X(iK+iB)=
(\pi\circ \mu)_*\cO_{X'}(iK_{X'}+iB'), \forall i\ge 1,
$$
$$
\bFix(iK+iB)=\bFix(iK_{X'}+iB')+i\overline{F}, \forall 
i\ge 1,
$$ 
and 
$a(E;Y,B_Y)=a(E;X',B'_Y)+\mult_E(\overline{F\vert_{Y'}})$
for every prime b-divisor $E$ of $Y$. The claim follows 
from the definition of $\bTheta$.
\end{proof}

\begin{lem}\label{2m} Let $\nu\colon Y''\to Y'$ be a 
birational contraction of birational models of $Y$. 
Then
$$
\nu_*\cO_{Y''}(nK_{Y''}+n\bTheta_{Y''})\subseteq
\cO_{Y'}(nK_{Y'}+n\bTheta_{Y'}), n\ge 1.
$$
Equality holds if $(Y',\bTheta_{Y'})$ is a log 
variety with canonical singularities in codimension 
at least two.
\end{lem}

\begin{proof} The inclusions are clear, since 
$K_{Y'}+\bTheta_{Y'}=\nu_*(K_{Y''}+\bTheta_{Y''})$.
For the second claim, there exists a $\nu$-exceptional 
$\R$-divisor $E$ such that
$$
K_{Y''}+\bTheta_{Y''}=\nu^*(K_{Y'}+\bTheta_{Y'})+E.
$$
If $\mult_P(E)<0$ for some prime divisor $P$ on 
$Y''$, then $a(P;Y',\bTheta_{Y'})<1$ since
$\bTheta_{Y''}$ is effective. Since $P$ is
$\nu$-exceptional, we infer that $(Y',\bTheta_{Y'})$ 
does not have canonical singularities in codimension at 
least two. Contradiction!
Therefore $E$ is effective and $\nu$-exceptional, which
implies the claim.
\end{proof}

\begin{lem}\label{af} Let $\nu\colon Y'\to Y$ be a 
birational contraction. Then 
\begin{itemize}
\item[(1)] For every $n\ge 1$, we have natural inclusions
$$
\im(\pi_*\cO_X(nK+nB)\to \pi_*\cO_Y(nK_Y+nB_Y))
\subseteq \pi_*\nu_*\cO_{Y'}(nK_{Y'}+n\bTheta_{Y'}).
$$
\item[(2)] The $\R$-divisor
$\lim_{i \to \infty}\frac{\Fix(iK_{Y'}+i\bTheta_{Y'})}{i}$ 
has zero multiplicity at each component of $\bTheta_{Y'}$.
\end{itemize}
\end{lem}

\begin{proof} By Lemma~\ref{2m}, we may assume that $Y'=Y$.
Let $i\in l\N$, so that $iK+iB$ is relatively mobile at $Y$.
By Hironaka, there exists a resolution $\mu_i\colon X_i\to X$
with the following properties:
\begin{itemize}
\item[(i)] Let $\mu^*(iK+iB)=M_i+F_i$ is the 
mobile-fixed decomposition relative to $S$. Then $M_i$ 
is relatively free.
\item[(ii)] Let $\mu^*(K+B)=K_{X_i}+B_{X_i}$ be the log
pullback. Then $\Supp(B_{X_i})\cup \Supp(F_i)$ is a simple
normal crossings divisor.
\end{itemize}
Let $Y_i$ be the proper transform of $Y$ on $X_i$.
Let $\pi_i=\pi\circ \mu$. We have 
$$
\pi_*\cO_X(nK+nB)={\pi_i}_*\cO_{X_i}(nK_{X_i}+nB_{X_i}).
$$

(1) Fix $n$ and choose $i\in ln\N$. We have 
$$
{\pi_i}_*\cO_{X_i}(nK_{X_i}+nB_{X_i}-\frac{n}{i}F_i)=
{\pi_i}_*\cO_{X_i}(nK_{X_i}+nB_{X_i}).
$$
Therefore the restriction to $Y_i$ of the right hand side 
is included in
$$
{\pi_i}_*\cO_{Y_i}(nK_{Y_i}+n(B_{X_i}-Y_i)\vert_{Y_i}-
\frac{n}{i}F_i\vert_{Y_i}).
$$
In turn, this sheaf is included in 
$$
\pi_*\cO_Y(nK_Y+nB_Y-(\mu\vert_{Y_i})_*
(\frac{n}{i}F_i\vert_{Y_i})).
$$
But $B_Y-(\mu\vert_{Y_i})_*(\frac{F_i}{i}\vert_{Y_i})\le 
\bTheta_Y$, hence the claim holds.

(2) For $i\in l\N$, denote 
$\tilde{B}_i=\max(B_Y-(\mu\vert_{Y_i})_*(\frac{F_i}{i}\vert_{Y_i}),0)$.
We claim that $iK_Y+i\tilde{B}_i$ is relatively mobile at 
the components of $\tilde{B}_i$.

Indeed, let $B_i=\max(B_{X_i}-\frac{F_i}{i},0)$. 
By construction, the mobile part of $iK_{X_i}+iB_i$ 
is $M_i$ and the fixed part $F'_i=iK_{X_i}+iB_i-M_i$ 
has no components in common with $B_i$. Since 
$\Supp(F'_i)\cup \Supp(B_i)$ is a simple normal crossings 
divisor, we infer that $F'_i$ does not contain any 
intersection of the components of $B_i$. In particular,
$iK_{Y_i}+i(B_i-Y_i)\vert_{Y_i}$ is relatively mobile at 
each component of $(B_i-Y_i)\vert_Y$. But 
$iK_Y+i\tilde{B}_i=\nu_*(iK_{Y_i}+i(B_i-Y_i)\vert_{Y_i})$, 
hence $iK_Y+i\tilde{B}_i$ is relatively mobile at the 
components of $\tilde{B}_i$.

By convexity, $\tilde{B}_i\le \bTheta_Y$ and 
$\lim_{i\to \infty}\tilde{B}_i=\bTheta_Y$.
Therefore the components of $\tilde{B}_i$ and $\bTheta_Y$ 
coincide for $i\gg 1$. We have 
$i(K_Y+\bTheta_Y)=i(K_Y+\tilde{B}_i)+i(\bTheta_Y-\tilde{B}_i)$.
Therefore $\Fix(iK_Y+i\bTheta_Y)\le i(\bTheta_Y-\tilde{B}_i)$ 
at each component of $\bTheta_Y$. Dividing 
by $i$ and taking the limit, we obtain the claim.
\end{proof}

\begin{lem}\label{cd} In the above notations, assume $(X,B)$ 
is a plt log variety with unique log canonical centre $Y$, and
$K+B\sim_\Q A+C$, where $A$ is a relatively ample 
$\R$-divisor and $C$ is an effective $\R$-divisor 
such that $\mult_Y(A)=\mult_Y(C)=0$.
\newline
If the (real) log Minimal Model Program holds in 
dimension $\dim(X)$, there exists a birational 
contraction $\nu\colon \tilde{Y}\to Y$ such that 
$$
\im(\pi_*\cO_X(nK+nB)\to \pi_*\cO_Y(nK_Y+nB_Y))=
\pi_*\nu_*\cO_{\tilde{Y}}(nK_{\tilde{Y}}+
n\bTheta_{\tilde{Y}})
$$
for every $n\ge 2$ such that $\{nB\}\le B$. Moreover,
$\bTheta$ is rational if $B$ is rational.
\end{lem}

\begin{proof} If the real logMMP holds, $(X,B)$ has a log
minimal model
$$
\Phi\colon (X,B)\dashrightarrow (X',B').
$$
If we consider a resolution of the graph of $\Phi$, this
means that we have a Hironaka hut
 \[ \xymatrix{
& \tilde{X} \ar[dl]_\mu \ar[dr]^{\mu'} &  \\
X \ar@{.>}[rr]^\Phi & & X'
 } \]
such that $K_{X'}+B'$ is relatively nef and big and the 
$\R$-Cartier divisor $F=\mu^*(K+B)-{\mu'}^*(K_{X'}+B')$ 
is effective and $\mu'$-exceptional. Denote $P={\mu'}^*(K_{X'}+B')$.
Then $\mu^*(K+B)=P+F$ is a Zariski decomposition:
\begin{itemize}
\item[(i)] $P\le \mu^*(K+B)$ is relatively nef and big;
\item[(ii)] $\cR_{\tilde{X}/S}(P)=\cR_{X/S}(K+B)$.
\end{itemize}
We may assume that $\tilde{Y}$, the proper transform of $Y$ 
on $\tilde{X}$, is normal. Since $lK+lB$ is relatively 
mobile at $Y$, we have $\mult_{\tilde{Y}}(F)=0$. We have
$$
P\vert_{\tilde{Y}}=K_{\tilde{Y}}+B_{\tilde{Y}}-F\vert_{\tilde{Y}}.
$$
By definition, 
$\bTheta_{\tilde{Y}}=\max(B_{\tilde{Y}}-F\vert_{\tilde{Y}},0)$. 
In particular,
$
P\vert_{\tilde{Y}}\le K_{\tilde{Y}}+\bTheta_{\tilde{Y}}.
$
We claim that this is a Zariski decomposition:
\begin{itemize}
\item[(iii)] $P\vert_{\tilde{Y}}\le K_{\tilde{Y}}+\bTheta_{\tilde{Y}}$
is relatively nef and big;
\item[(iv)] $\cR_{\tilde{Y}/S}(P\vert_{\tilde{Y}})=\cR_{\tilde{Y}/S}
(K_{\tilde{Y}}+\bTheta_{\tilde{Y}})$.
\end{itemize}
Indeed, by assumption $Y$ is mapped birationally by $\Phi$ 
to a prime divisor $Y'$. Since $(X',B')$ is also plt, $Y'$
is normal. Let $(Y',B'_{Y'})$ be the log variety structure 
induced by codimension one adjunction. By assumption, 
$K_{Y'}+B'_{Y'}$ is relatively nef and big. By adjunction, 
we have an induced Hironaka hut
 \[ \xymatrix{
& \tilde{Y} \ar[dl]_\nu \ar[dr]^{\nu'} &  \\
Y \ar@{.>}[rr]^{\Phi\vert_Y} & & Y'.
 } \]
We have ${\nu'}^*(K_{Y'}+B'_{Y'})=K_{\tilde{Y}}+B_{\tilde{Y}}-
F\vert_{\tilde{Y}}$. Since $B'_{Y'}$ is effective, the
negative part of $B_{\tilde{Y}}-F\vert_{\tilde{Y}}$ is 
$\nu'$-exceptional. Therefore 
$
P\vert_{\tilde{Y}}\le K_{\tilde{Y}}+\bTheta_{\tilde{Y}}.
$
is a Zariski decomposition.

Let $n\ge 2$ be an integer and let $i\in nl\N$. 
We have a commutative diagram
\[ \xymatrix{
\pi_*\cO_X(nK+nB)   \ar[d]  \ar[r] & 
\pi_*\cO_X(\lceil \bK+Y+(n-1)\bD_i(K+B)\rceil)   \ar[d]\\
\pi_*\cO_{\tilde{Y}}(nK_{\tilde{Y}}+n\bTheta_{\tilde{Y}})  \ar[r]    & 
\pi_*\cO_Y(\lceil \bK+(n-1)
\bD_i(K_{\tilde{Y}}+\bTheta_{\tilde{Y}}))\rceil)
} \]
The horizontal arrows are inclusions. For $i$ sufficiently
large and divisible, the right hand side vertical arrow becomes
$$
\pi_*\cO_X(\lceil \bK+Y+(n-1)\overline{P}\rceil)
\to
\pi_*\cO_{\tilde{Y}}(\lceil \bK+(n-1)\overline{P\vert_{\tilde{Y}}}\rceil),
$$
which is surjective by Kawamata-Viehweg vanishing. Assume 
moreover that $\{nB\}\le B$. Then the top horizontal arrow 
is an equality, by Lemma~\ref{k+}. Therefore the 
bottom horizontal arrow is an equality and the left hand 
side vertical arrow is surjective. 

Finally, assume that $B$ is rational. Then $B'$ is rational, 
hence $F$ is rational. Therefore $\bTheta$ is rational.
\end{proof}

%%%%%%%%%%%%%%%%%%%%%%%%%%%%%%%%%%%%%%%
%%%%%%%%%%%%%%%%%%%%%%%%%%%%%%%%%%%%%%%

\section{Proof of Theorem~\ref{mr}}

%%%%%%%%%%%%%%%%%%%%%%%%%%%%%%%%%%%%%%%
%%%%%%%%%%%%%%%%%%%%%%%%%%%%%%%%%%%%%%%

First of all, by Lemmas~\ref{cut} and~\ref{2m}, we may assume 
that the prime components of $B_Y$ are disjoint. Note that 
assumption (a) is birational in nature (see~\cite{Nak98}, 
Lemma 4.8). Denote $\Theta=\bTheta_Y$.

Let $H$ be a $\pi$-very ample divisor with $\mult_Y(H)=0$
and let $A=\dim(X)\cdot H$. There exists a positive integer 
$l$ such that $l(K+B)\sim A+C$, where $C$ is an 
effective $\R$-divisor with $\mult_Y(C)=0$. 
For $n\ge 1$, define
$$
\begin{array}{l}
B_n=\max(B,\frac{1}{n}\lceil (n-1)B \rceil)\\
B_{n,Y}=(B_n-Y)\vert_Y\\
\Theta_n=\max(B_{n,Y}-\lim_{i\to \infty}
\frac{(\bFix(i(K+B_n+\frac{1}{n}A))\vert_Y)_Y}{i},0)\\
C_n=C+l(B_n-B)
\end{array}
$$
Note that $(X,B_n)$ satisfies the same properties 
as $(X,B)$.

\begin{lem}\label{am} The following properties hold:
\begin{itemize}
\item[(1)] $Y+\lceil (n-1)B_{n-1}\rceil\le n B_n$ 
for $n\ge 1$.
\item[(2)] $\Theta\le \Theta_n\le B_{n,Y}$.
\item[(3)] $\bD_i(nK+nB_n+A)\vert_Y\le 
\bD_i(nK_Y+n\Theta_n+A\vert_Y)$ for $i\in l\N$. 
\end{itemize}
\end{lem}

\begin{proof} (1) This follows from Lemma~\ref{Sep22}.

(2) The inequality $\Theta_n\le B_{n,Y}$ holds 
by construction. We have
$$
K+B_n+\frac{1}{n}A=(K+B)+(B_n-B)+\frac{1}{n}A.
$$
Since $B_n-B$ is effective and $A$ is $\pi$-free, 
we obtain
$$
\bFix(i(K+B_n+\frac{1}{n}A))\le \bFix(i(K+B))
+i\overline{(B_n-B)}\ \forall i\in ln\N,
$$
which is equivalent to
$$
\overline{B-Y}-\frac{\bFix(i(K+B))}{i}
\le \overline{B_n-Y}-\frac{\bFix(i(K+B_n+\frac{1}{n}A))}{i}.
$$
Restricting these $\R$-Cartier $\R$-b-divisors to $Y$ and 
taking the trace on $Y$ we obtain
$$
B_Y-\frac{\bFix(i(K+B))}{i}\vert_Y
\le B_{n,Y}-\frac{\bFix(i(K+B_n+\frac{1}{n}A))}{i}\vert_Y.
$$
Taking the effective part and the limit with respect to $i$, 
we obtain $\Theta\le \Theta_n$. 

(3) Let $i\in ln\N$. By definition, the mobile b-divisor of
$i(K+B_n+\frac{1}{n}A)$ coincides with the mobile b-divisor 
of
$$
i\overline{(K+B_n+\frac{1}{n}A)}-
\bFix(i(K+B_n+\frac{1}{n}A)).
$$
The restriction of this $\R$-Cartier $\R$-b-divisor to $Y$ is 
$$
i\overline{(K_Y+B_{n,Y}+\frac{1}{n}A\vert_Y)}-
\bFix(i(K+B_n+\frac{1}{n}A))\vert_Y,
$$
whose mobile b-divisor is at most the mobile b-divisor of 
$$
i(K_Y+B_{n,Y}+\frac{1}{n}A\vert_Y)-
(\bFix(i(K+B_n+\frac{1}{n}A))\vert_Y)_Y.
$$
By the definition of $\Theta_n$ and the convexity of the 
sequence of fixed parts, we have 
$$
B_{n,Y}-\frac{(\bFix(i(K+B_n+\frac{1}{n}A))\vert_Y)_Y}{i}
\le \Theta_n.
$$
Therefore the restriction to $Y$ of the mobile b-divisor of
$i(K+B_n+\frac{1}{n}A)$ is at most the mobile b-divisor
of $i(K_Y+\Theta_n)$. We obtain the claim by dividing by $i$.
\end{proof}

\begin{proof}[Proof of Theorem \eqref{mr}$(1)$]  
{\em Step 1.} For every $n\ge 0$, there exists $i_n\in l\N$ such
that for every $i\in i_n\N$ the following inclusion holds
$$
\cO_Y(\lceil \bK+n\bD_i(K_Y+\Theta)+\overline{A\vert_Y}\rceil)
\subseteq
\cO_Y(\lceil \bK+\bD_i(n(K+B_n)+A)\vert_Y\rceil).
$$
For $n=0$, set $i_0=l$. Since $A$ is free, we have 
$\bD_i(A)=\overline{A}$ for every $i$. Therefore the inclusion 
is an equality.

Let now $n\ge 1$ and assume the inclusion holds for $n-1$, with
corresponding index $i_{n-1}$. Let $i_n$ be the smallest positive
integer $j$ satisfying the following properties:
\begin{itemize}
\item[(i)] $i_{n-1}\vert j$,
\item[(ii)] $\bD_j((n-1)(K+B_{n-1})+A)$ is b-big,
\item[(iii)] $\lceil n(K_Y+\Theta_n)+A\vert_Y\rceil\le
\lceil \bD_j(n(K+B_n)+A)\vert_Y\rceil$ at every prime
component of $\Theta$.
\end{itemize}
The existence of a solution for (iii) follows from the 
inequality $\Theta\le \Theta_n$ and the definition of 
$\Theta_n$.
Let $i\in i_n\N$. Let $R$ be the reduced divisor supported by 
$B-Y$ with $R\vert_Y=\lceil \Theta\rceil$. The following 
inclusions hold:
\begin{align*}
& \pi_*\cO_X(\lceil \bK+Y+\bD_i((n-1)(K+B_{n-1})+
\overline{A})\rceil) \\
&\subseteq  \pi_*\cO_X(\lceil K+Y+(n-1)(K+B_{n-1})+A\rceil) \\
&\subseteq  \pi_*\cO_X(n(K+B_n)+A)  \\
&\subseteq  \pi_*\cO_X(\lceil \bA(X,Y+R)^*+\bD_i(n(K+B_n)+A)\rceil)
\end{align*}
Indeed, the first inclusion is clear, the second holds by
Lemma~\ref{am}.(1), and the third follows from the inequality
$\lceil \bA(X,Y+R)^*\rceil\ge 0$. By Kawamata-Viehweg vanishing,
the restriction to $Y$ of the first sheaf is 
$\pi_*\cO_Y(\lceil \bK+\bD_i((n-1)(K+B_{n-1})+A)\vert_Y\rceil)$.
Since $\bA(X,Y+R)^*\vert_Y=\bA(Y,\lceil \Theta\rceil)^*$,
we obtain by restricting to $Y$
\begin{align*}
& \pi_*\cO_Y(\lceil \bK+\bD_i((n-1)(K+B_{n-1})+A)\vert_Y\rceil) \\
& \subseteq 
\pi_*\cO_Y(\lceil \bA(Y,\lceil \Theta\rceil)^*+\bD_i(n(K+B_n)+A)
\vert_Y\rceil).
\end{align*}
Combining this with Step 1 for $n-1$, we obtain
\begin{align*}
& \pi_*\cO_Y(\lceil \bK+(n-1)\bD_i(K_Y+\Theta)+\overline{A\vert_Y}\rceil) \\
& \subseteq 
\pi_*\cO_Y(\lceil \bA(Y,\lceil \Theta\rceil)^*+\bD_i(n(K+B_n)+A)
\vert_Y\rceil).
\end{align*}
Since $\cO_Y(\lceil \bK+(n-1)\bD_i(K_Y+\Theta)+\overline{A\vert_Y}\rceil)$
is $\pi\vert_Y$-generated (\cite{Nak98}, Lemma 3.9), this is equivalent to 
\begin{align*}
& \cO_Y(\lceil \bK+(n-1)\bD_i(K_Y+\Theta)+\overline{A\vert_Y}\rceil) \\
& \subseteq 
\cO_Y(\lceil \bA(Y,\lceil \Theta\rceil)^*+\bD_i(n(K+B_n)+A)
\vert_Y\rceil).
\end{align*}
In particular, we obtain inclusions
\begin{align*}
&\cO_Y(\lceil \bK+n\bD_i(K_Y+\Theta)+\overline{A\vert_Y}\rceil) \\
&\subseteq \cO_Y(\lceil \bK+(n-1)\bD_i(K_Y+\Theta)+\overline{A\vert_Y}\rceil
+\overline{K_Y+\lceil \Theta\rceil})\\
&\subseteq 
\cO_Y(\lceil \bA(Y,\lceil \Theta\rceil)^*+
\overline{K_Y+\lceil \Theta\rceil}+
\bD_i(n(K+B_n)+A)\vert_Y \rceil)\\
&=\cO_Y(\lceil \bK+\lceil \Theta\rceil+\bD_i(n(K+B_n)+A)\vert_Y \rceil).
\end{align*}
Let $U\subset Y$ be an open set and let $a$ be a nonzero
rational function on $Y$ such that 
$$
\overline{(a)}+\lceil\bK+n\bD_i(K_Y+\Theta)+\overline{A\vert_Y}
\rceil\vert_U\ge 0. 
$$
We can write
$
\overline{(a)}+\lceil\bK+\bD_i(n(K+B_n)+A)\vert_Y
\rceil\vert_U=\bE^+-\bE^-,
$
where $\bE^+,\bE^-$ are effective b-divisors of $U$ with
no common components. From the inclusions above, we have
$\lceil \Theta\rceil+\bE^+-\bE^-\ge 0$. Therefore
$\bE^-\le \lceil \Theta\rceil$. But $\bE^-$ has zero multiplicity
at the components of $\Theta$, by (iii). Therefore $\bE^-=0$, that is
$$
\overline{(a)}+\lceil\bK+\bD_i(n(K+B_n)+A)\vert_Y
\rceil\vert_U\ge 0.
$$
This shows that 
$$
\cO_Y(\lceil \bK+n\bD(K_Y+\Theta)+\overline{A\vert_Y}\rceil)
\subseteq
\cO_Y(\lceil \bK+\bD(n(K+B_n)+A)\vert_Y\rceil),
$$
which is the desired inclusion for $n$.

{\em Step 2.} For $n\ge 1$ and $i\in i_n\N$ we have
$$
\cO_Y(\lceil \bK+(n-1)\bD_i(K_Y+\Theta)+\overline{A\vert_Y}\rceil)
\subseteq
\cO_Y(\lceil \bA(Y,B_{n,Y})+\bD_i(n(K+B_n)+A)\vert_Y\rceil).
$$
To see this, note first the inclusions
\begin{align*}
&\pi_*\cO_X(\lceil \bK+Y+\bD_i((n-1)(K+B_{n-1})+\overline{A})\rceil) \\
&\subseteq  \pi_*\cO_X(\lceil K+Y+(n-1)(K+B_{n-1})+A\rceil) \\
&\subseteq  \pi_*\cO_X(n(K+B_n)+A)  \\
&\subseteq  \pi_*\cO_X(\lceil \bA(X,B_n)^*+\bD_i(n(K+B_n)+A)\rceil).
\end{align*}
The argument is the same as in Step 1, except that for the
last inclusion we use $\lceil \bA(X,B_n)^*\rceil\ge 0$.
By Kawamata-Viehweg vanishing, the restriction to $Y$ of the
first sheaf is $\pi_*\cO_Y(\lceil \bK+(n-1)\bD_i(K_Y+\Theta)+
\overline{A\vert_Y}\rceil)$.
Since $\bA(X,B_n)^*\vert_Y=\bA(Y,B_{n,Y})^*$, we obtain
\begin{align*}
& \pi_*\cO_Y(\lceil \bK+\bD_i((n-1)(K+B_{n-1})+A)\vert_Y\rceil) \\
& \subseteq 
\pi_*\cO_Y(\lceil \bA(Y,B_{n,Y})+\bD_i(n(K+B_n)+A)\vert_Y\rceil)
\end{align*}
Combining this with Step 1, we obtain
\begin{align*}
& \pi_*\cO_Y(\lceil \bK+(n-1)\bD_i(K_Y+\Theta)+\overline{A\vert_Y}\rceil) \\
& \subseteq 
\pi_*\cO_Y(\lceil \bA(Y,B_{n,Y})+\bD_i(n(K+B_n)+A)\vert_Y\rceil)
\end{align*}
Since $\cO_Y(\lceil \bK+(n-1)\bD_i(K_Y+\Theta)+\overline{A\vert_Y}\rceil)$ 
is $\pi\vert_Y$-generated (\cite{Nak98}, Lemma 3.9), 
this is equivalent to the claim.

{\em Step 3.} For $n\ge 1$ and $i\in i_n\N$ we have
$$
\cO_Y(\lceil \bK+(n-1)\bD_i(K_Y+\Theta)\rceil)
\subseteq \cO_Y(\lceil \bA(Y,B_{n,Y})+
\overline{C_n\vert_Y}+n\bD_i(K+B_n)\vert_Y\rceil).
$$
Indeed, $l(K+B_n)\sim A+C_n$. Since $l\ge 1$ and $A$ is 
$\pi$-free, we obtain
$$
\bD_i(n(K+B_n)+A)\le n\bD_i(K+B_n)+\overline{A+C_n}.
$$
We obtain the claim by restricting to $Y$, using Step 2 
and canceling $A\vert_Y$. 

{\em Step 4.} For $n\ge 1$, there exists $i'_n\in i_n\N$
such that for every $i\in i'_n\N$ we have
$$
\pi_*\cO_Y(n(K_Y+\Theta))\subseteq \pi_*\cO_Y(\lceil 
\bA(Y,B_{n,Y})+n\bD_i(K+B_n)\vert_Y\rceil).
$$
Fix $n\ge 1$. There exists $\gamma_n>0$ such that the log variety
$$
(Y,B_{n,Y}+\gamma(C\vert_Y+(\frac{l}{n}+1)\lceil B_Y\rceil))
$$ 
has Kawamata log terminal singularities. By diophantine 
approximation~\cite{Cassels57},
there exists an integer $e\ge 1+\gamma^{-1}$ such that 
$B_{ne}\le B_n$. Define $i'_n=i_{ne}$. 
Let $\bM_n$ be the mobile b-divisor of $\pi_*\cO_Y(nK_Y+n\Theta)$. 
By Step 3, for every $i\in i'_n$ we have 
$$
\bM_n\le \frac{1}{e}\bM_{ne}\le
\frac{1}{e}\lceil \bA(Y,B_{ne,Y})+\overline{C_n\vert_Y}+
ne\bD_i(K+B_{ne})\vert_Y \rceil. 
$$
Since $\bM_n$ has integer coefficients, this is equivalent to
$$
\bM_n\le \lfloor  
\frac{1}{e}\lceil \bA(Y,B_{n,Y})+
\overline{(B_n-B_{ne}+C_n)\vert_Y}+
ne\bD_i(K+B_{ne})\vert_Y\rceil
\rfloor
$$
By the choice of $\gamma$ and Lemma~\ref{Sep22}.(2), we have 
$$
\lceil \bA(Y,B_{n,Y})-\gamma
\overline{(B_n-B_{ne}+C_n)\vert_Y} \rceil \ge 0.
$$
Therefore we may apply Lemma~\ref{div} at each prime b-divisor of
$Y$ and obtain
$$
\bM_n\le
\lceil \bA(Y,B_{n,Y})+n\bD_i(K+B_{ne})\vert_Y \rceil. 
$$
We have $\bD_i(K+B_{ne})\le \bD_i(K+B_n)$ since $B_{ne}\le B_n$, 
hence we get the claim.

{\em Step 5.} For $n\ge 2$ we have
$$
\pi_*\cO_Y(nK_Y+n\Theta)\subseteq
\im(\pi_*\cO_X(nK+nB_n)\to \pi_*\cO_Y(nK_Y+nB_{n,Y})).
$$
Indeed, the following inequality holds:
$$
\bA(Y,B_{n,Y})+n\bD_i(K+B_n)\vert_Y\le
\bK+(n-1)\bD_i(K+B_n)\vert_Y.
$$
By Step 4, for $n\ge 1$ and $i\gg 1$ we have
$$
\pi_*\cO_Y(nK_Y+n\Theta)\subseteq 
\pi_*\cO_Y(\lceil \bK+(n-1)\bD_i(K+B_n)\vert_Y\rceil).
$$
By Kawamata-Viehweg vanishing, 
$\pi_*\cO_Y(\lceil \bK+(n-1)\bD_i(K+B_n)\vert_Y)$ lifts to
$
\pi_*\cO_X(\lceil\bK+Y+(n-1)\bD_i(K+B_{n-1})\rceil),
$
which is included in $\pi_*\cO_X(\lceil K+Y+(n-1)(K+B_{n-1})\rceil)$.
By Lemma~\ref{Sep22}.(3), this is included in 
$\pi_*\cO_X(nK+nB_n)$. This proves the claim.

{\em Step 6.} Let $n\ge 2$ such that $\{nB\}\le B$. Then
$$
\pi_*\cO_Y(nK_Y+n\Theta)=
\im(\pi_*\cO_X(nK+nB)\to \pi_*\cO_Y(nK_Y+nB_Y)).
$$
The second assumption means $B_n=B$, hence
Step 5 gives the direct inclusion. The opposite 
inclusion is clear.

{\em Step 7.} Assume that $n=1$ and $\pi(Y)\ne \pi(X)$.
Then 
$$
\im(\pi_*\cO_X(K+B)\to \pi_*\cO_Y(K_Y+B_Y))=\pi_*\cO_Y(K_Y+\Theta).
$$
Indeed, Steps 5 and 6 hold for $n=1$ as well, by Koll\'ar's
torsion freeness rather than Kawamata-Viehweg vanishing.
\end{proof}

\begin{prop}[Generalized asymptotic saturation]\label{sat}
Let $n\ge 2$ such that $\{nB\}\le B$ and let $i\ge 1$ such that
$\pi_*\cO_X(iK+iB)\ne 0$. Then
$$
\pi_*\cO_Y(\lceil \bK+(n-1)\bD_i(K_Y+\Theta)\rceil)
\subseteq \pi_*\cO_Y(nK_Y+n\Theta).
$$
\end{prop}

\begin{proof} The log canonical divisor $K_Y+\Theta$ is
relatively big, by the assumption (a). By diophantine 
approximation~\cite{Cassels57}, there exists an integer
$i'\ge 2$ such that $i\vert i'$, $\{i'B\}\le B$ and 
$\bD_{i'}(K+B)$ is relatively b-nef and b-big. 
By Theorem~\ref{mr}.(1), we have
$
\bD_{i'}(K+B)\vert_Y=\bD_{i'}(K_Y+\Theta).
$
In particular, 
$$
\pi_*\cO_Y(\lceil \bK+(n-1)\bD_{i'}(K_Y+\Theta)\rceil)=
\pi_*\cO_Y(\lceil \bK+(n-1)\bD_{i'}(K+B)\vert_Y\rceil).
$$
Since $\bD_{i'}(K+B)$ is relatively b-nef and b-big, 
the right hand side lifts to 
$\pi_*\cO_Y(\lceil \bK+Y+(n-1)\bD_{i'}(K+B)\rceil)$,
by Kawamata-Viehweg vanishing. Since $\{nB\}\le B$, we have
$$
\pi_*\cO_Y(\lceil \bK+Y+(n-1)\bD_i(K+B)\rceil)
\subseteq 
\pi_*\cO_X(nK+nB).
$$
The restriction to $Y$ of $\pi_*\cO_X(nK+nB)$ is 
$\pi_*\cO_Y(nK_Y+n\Theta)$, by Theorem~\ref{mr}.(1).
Therefore
$$
\pi_*\cO_Y(\lceil \bK+(n-1)\bD_{i'}(K_Y+\Theta)\rceil)
\subseteq \pi_*\cO_Y(nK_Y+n\Theta).
$$
The claim follows from
$\bD_i(K_Y+\Theta)\le \bD_{i'}(K_Y+\Theta)$.
\end{proof}

\begin{proof}[Proof of Theorem \eqref{mr}$(2)$]  
Let $I\ge 2$ be an integer such that $IB$ is integral. 
Assume that the real log divisor $K_Y+\Theta$, 
which is $\pi\vert_Y$-big by the assumption (a) in 
Theorem~\ref{mr}, has a Zariski decomposition. Thus 
there exists a birational contraction $\mu\colon Y'\to Y$ 
and a $(\pi\circ\mu)$-nef and $(\pi\circ\mu)$-big 
$\R$-divisor $P$ on $Y'$ such that 
\begin{itemize}
\item[(i)] $P\le \mu^*(K_Y+\Theta)$;
\item[(ii)] $\cR_{Y'/S}(P)=\cR_{Y/S}(K_Y+\Theta)$.
\end{itemize}

{\em Step 1.} We claim that $P$ is rational and 
relatively semiample. In particular, by (ii), the 
$\cO_S$-algebra $\cR_{Y/S}(K_Y+\Theta)$ is finitely 
generated.

Indeed, after possibly blowing up $Y'$ and replacing $P$ 
by its pullback, we may assume that $Y'$ is smooth and 
$\Supp(P)\cup \Supp(\Theta_{Y'})$ is a simple normal 
crossings divisor, where $\mu^*(K_Y+\Theta)=K_{Y'}+\Theta_{Y'}$ 
is the log pullback. Denote $\pi'=\pi\circ \mu$.
By (ii) and the argument of~\cite{Semcr}, Proposition 3.1, we 
have $\lim_{i\to \infty}\bD_i(K_Y+\Theta)=\overline{P}$ and the 
inclusions in Lemma~\ref{sat} become
$$
\pi_*\cO_Y(\lceil \bK+(n-1)\overline{P}\rceil)\subseteq
\pi_*\cO_Y(nK_Y+n\Theta), \forall I\vert n.
$$
Since $\pi'_*\cO_{Y'}(P)=\pi_*\cO_Y(nK_Y+n\Theta)$ and 
$\Supp(P)$ has simple normal crossings, this is equivalent
to
$$
\pi'_*\cO_{Y'}(\lceil K_{Y'}+(n-1)P\rceil)\subseteq \pi'_*\cO_{Y'}(nP), 
\forall I\vert n.
$$
If we denote $N=K_{Y'}+\Theta_{Y'}-P\ge 0$, these inclusions
become
$$
\pi'_*\cO_{Y'}(\lceil -(\Theta_{Y'}-N)+n P\rceil)\subseteq 
\pi'_*\cO_{Y'}(nP), \forall I\vert n.
$$ 
The log pair $(Y',\Theta_{Y'}-N)$ has Kawamata log terminal 
singularities (the boundary may not be effective), and
$2P-(K_{Y'}+\Theta_{Y'}-N)=P$ is a $\pi'$-nef and $\pi'$-big
$\R$-divisor. We have verified the assumptions of~\cite{Semcr}, 
Theorem 2.1 for $(Y'/S, \Theta_{Y'}-N)$ and $P$, hence we conclude
that $P$ is rational and $\pi'$-semiample.

{\em Step 2.} Since $P$ is relatively semiample and big, there 
exists a birational contraction $\nu\colon Y'\to Y''$, defined 
over $S$, and a relatively ample $\Q$-Cartier divisor $P''$ on 
$Y''$ such that $P=\nu^*(P'')$.
\[ \xymatrix{
& Y' \ar[dl]_\mu \ar[dr]^{\nu} &  \\
Y \ar@{.>}[rr] & & Y''
 } \]
Let $\pi''\colon Y''\to S$ be the induced morphism. Since
$P''$ is relatively ample, there exists a positive integer
$n$ such that $nP''$ is Cartier and the sheaf
$\nu_*\cO_{Y'}(nN)\otimes\cO_{Y''}(nP'')$ is $\pi''$-generated.
We have $\pi'_*\cO_{Y'}(nP)=\pi'_*\cO_{Y'}(nP+nN)$, hence 
$$
\pi''_*\cO_{Y'}(nP'')=
\pi''_*(\nu_*\cO_{Y'}(nN)\otimes\cO_Y''(nP'')).
$$
The right hand side is generated by global sections,
hence 
$$
\nu_*\cO_{Y'}(nN)\otimes\cO_{Y''}(nP'')\subseteq
\cO_{Y'}(nP''),
$$
that is $\nu_*\cO_{Y'}(nN)\subseteq \cO_{Y''}$. Therefore
$N$ is $\nu$-exceptional. In particular, 
$$
P''=K_{Y''}+\nu_*(\Theta_{Y'}),
$$
and $\nu_*(\Theta_{Y'})$ is rational. We have
$$
\mu^*(K_Y+\Theta)=\nu^*(K_{Y''}+\nu_*(\Theta_{Y'}))+N.
$$
In particular, we have $\Fix(\mu^*(iK_Y+i\Theta))\ge iN$ for
every $i$, hence 
$\lim_{i\to \infty} \frac{\Fix(\mu^*(iK_Y+i\Theta))}{i}\ge N$.
Pushing this down to $Y$, we obtain
$$
\lim_{i\to \infty} \frac{\Fix(iK_Y+i\Theta)}{i}\ge \mu_*N.
$$
By Lemma~\ref{af}.(2), $\mu_*N$ has zero multiplicity 
at each component of $\Theta$. Combining with the above, 
we obtain
$$
\mult_E(\Theta)=\mult_E(\nu^*(K_{Y''}+\nu_*(\Theta_{Y'}))-K_{Y'})
$$
for each prime divisor $E$ in the support of $\Theta$. 
Therefore $\Theta$ is rational.
\end{proof}

%%%%%%%%%%%%%%%%%%%%%%%%%%%%%%%%%%%%%%%
%%%%%%%%%%%%%%%%%%%%%%%%%%%%%%%%%%%%%%%

\section{Applications}

%%%%%%%%%%%%%%%%%%%%%%%%%%%%%%%%%%%%%%%
%%%%%%%%%%%%%%%%%%%%%%%%%%%%%%%%%%%%%%%

Theorem~\ref{mth} generalizes~\cite{JH05}, Theorem 4.3,
and simplifies its proof. In particular, Shokurov's
prelimiting flips~\cite{Sh03} exist if big (real) log 
canonical divisors have a Zariski decomposition in one 
dimension less. For the convenience of the reader,
we compare the two approaches. With different assumptions,
Hacon and M\textsuperscript{c}Kernan prove Theorem~\ref{mth}.(1) 
in two steps. The first step proves (1) in the special case 
$\Theta=B_Y$. In particular,
the restricted algebra $\cR_{X/S}(K+B)\vert_Y$ has a 
representation
$\bigoplus_{n=0}^\infty \pi_*\cO_Y(nK_Y+n\Delta_n)$, where
$(\Delta_n)_{n\ge 0}$ is a sequence of boundaries
converging to $\Theta$. The second step shows the stabilization
$\Delta_n=\Theta$ for suitable $n$, and here one needs the log
Minimal Model Program in the dimension of $Y$, Shokurov's
asymptotic saturation of the restricted algebra, and the
log Fano assumption that $-(K+B)$ is $\pi$-ample.
The log Fano assumption is essential in this second step,
in order to apply Shokurov's asymptotic saturation, and it 
is unclear how to proceed without it. Instead, using the
new diophantine ingredients in \S 1,
we obtain (1) directly. As for Theorem~\ref{mth}.(2), Shokurov's 
asymptotic saturation is not enough for rationality
and finite generation if $-(K+B)$ is not $\pi$-ample. 
Instead, we introduce a stronger 
diophantine property (Lemma~\ref{k+}), which 
is preserved by restriction to $Y$ (Proposition~\ref{sat}). 

\begin{thm}\label{mth} 
Let $(X,B)$ be a log variety endowed with a proper 
contraction $\pi\colon X\to S$, satisfying the 
following properties:
\begin{itemize}
\item[(a)] $(X,B)$ is pure log terminal, with 
unique log canonical center $Y\subset X$.
\item[(b)] $K+B\sim_\Q A+C$, where $A$ is a $\pi$-ample
$\R$-divisor and $C$ is an effective $\R$-divisor with 
$\mult_Y(A)=\mult_Y(C)=0$.
\end{itemize}
Thus $Y$ is normal, and by codimension one adjunction
there exists a canonically defined boundary $B_Y$ such 
that 
$
(K+B)\vert_Y=K_Y+B_Y.
$
Then there exists a birational modification 
$\mu\colon \tilde{Y}\to Y$ with a natural structure 
$(\tilde{Y},\Theta)$ of log variety with Kawamata log 
terminal singularities, satisfying the following properties:
\begin{itemize}
\item[(1)] For every $n\ge 1$ we have natural inclusions
$$
\im(\pi_*\cO_X(nK+nB)\to \pi_*\cO_Y(nK_Y+nB_Y))\subseteq
(\pi\circ\mu)_*\cO_{\tilde{Y}}(nK_{\tilde{Y}}+n\Theta).
$$
The inclusion is an equality if $n\ge 2$ and $\{nB\}\le B$.
\item[(2)] Assume that $B$ has rational coefficients and
the log canonical divisor $K_{\tilde{Y}}+\Theta$ has a Zariski 
decomposition relative to $S$. 
Then $\Theta$ has rational coefficients and 
$\bigoplus_{n=0}^\infty (\pi\circ\mu)_*\cO_{\tilde{Y}}
(nK_{\tilde{Y}}+n\Theta)$ is a finitely generated
$\cO_S$-algebra.
\end{itemize}
\end{thm}

\begin{proof} The log variety $(Y,B_Y)$ has Kawamata
log terminal singularities. By Lemma~\ref{fn}, it
has finitely many valuations with log discrepancy
less than one. By Hironaka's resolution of singularities, 
we may construct a log resolution 
\[ \xymatrix{
\tilde{Y} \ar[r]\ar[d]_\nu    & \tilde{X} \ar[d]^\mu \\
Y    \ar[r]    &   X
} \] 
with the following properties:
\begin{itemize}
\item[(1)] The exceptional locus $\Exc(\mu)$ is a divisor
on $\tilde{X}$ and $\Exc(\mu)\cup \Supp(B_{\tilde{X}})$ has
simple normal crossings, where 
$\mu^*(K+B)=K_{\tilde{X}}+B_{\tilde{X}}$ is the log pullback. 
\item[(2)] $\tilde{Y}$ is the proper transform of
$Y$ on $\tilde{X}$, with $\mult_{\tilde{Y}}(B_{\tilde{X}})=1$.
\item[(3)] Let $\tilde{B}=\max(B_{\tilde{X}},0)$ and
$\tilde{B}_{\tilde{Y}}=(\tilde{B}-\tilde{Y})\vert_{\tilde{Y}}$. 
Then the components of $\tilde{B}_{\tilde{Y}}$ are disjoint.
\end{itemize}
Define
\begin{align*}
\Theta & = \max(B_{\tilde{Y}}-\lim_{i\to \infty}
\frac{(\bFix(iK+iB)\vert_{\tilde{Y}})_{\tilde{Y}}}{i},0) \\
  & =\max(\tilde{B}_{\tilde{Y}}-\lim_{i\to \infty}
\frac{(\bFix(iK_{\tilde{X}}+i\tilde{B})
\vert_{\tilde{Y}})_{\tilde{Y}}}{i},0),
\end{align*}
where the second equality holds by Lemma~\ref{cut}.
Let $\tilde{\pi}=\pi\circ \mu$. By Theorem~\ref{mr}.(1), 
we have inclusions
$$
\im(\tilde{\pi}_*\cO_{\tilde{X}}(nK_{\tilde{X}}+n\tilde{B})
\to \tilde{\pi}_*\cO_{\tilde{Y}}(nK_{\tilde{Y}}+n\tilde{B}
_{\tilde{Y}}))
\subseteq 
\tilde{\pi}_*\cO_{\tilde{Y}}(nK_{\tilde{Y}}+n\Theta)
$$
for every $n\ge 1$, and equality holds if $n\ge 2$ and
$\{n\tilde{B}\}\le \tilde{B}$.

Let now $n\ge 2$ with $\{nB\}\le B$. In particular, 
$$
\pi_*\cO_X(nK+nB)=\pi_*\cO_X(\lceil K+Y+(n-1)(K+B)\rceil).
$$
Therefore
$$
\pi_*\cO_X(nK+nB)=\pi_*\cO_X(\lceil \bK+Y+(n-1)\bD_i(K+B)\rceil),
i\gg 1.
$$
This is equivalent to
$$
\tilde{\pi}_*\cO_{\tilde{X}}(nK_{\tilde{X}}+n\tilde{B})
=\tilde{\pi}_*\cO_{\tilde{X}}(\lceil \bK+\tilde{Y}+
(n-1)\bD_i(K_{\tilde{X}}+\tilde{B})\rceil),
i\gg 1.
$$
We have a commutative diagram
\[ \xymatrix{
\tilde{\pi}_*\cO_{\tilde{X}}(nK_{\tilde{X}}+n\tilde{B})   \ar[d]  \ar[r] & 
\tilde{\pi}_*\cO_{\tilde{X}}(\lceil \bK+\tilde{Y}+
(n-1)\bD_i(K_{\tilde{X}}+\tilde{B})\rceil)   \ar[d]\\
\tilde{\pi}_*\cO_{\tilde{Y}}(nK_{\tilde{Y}}+n\Theta)  \ar[r]    & 
\tilde{\pi}_*\cO_{\tilde{Y}}(\lceil \bK+(n-1)
\bD_i(K_{\tilde{Y}}+\Theta))\rceil)
} \]
The horizontal arrows are inclusions. Choose
$i\ge 2$ such that $n\vert i$, $\{i\tilde{B}\}\le \tilde{B}$
and $\bD_i(K_{\tilde{X}}+\tilde{B})$ is relatively b-big.
From above, we have
$
\bD_i(K_{\tilde{X}}+\tilde{B})\vert_{\tilde{Y}}=
\bD_i(K_{\tilde{Y}}+\Theta),
$
hence by Kawamata-Viehweg vanishing we infer that
the right hand side vertical arrow is surjective.
The top horizontal arrow is an equality from above. This 
implies that the bottom horizontal arrow is an equality 
and the left hand side vertical arrow is surjective. 
This finishes the proof of (1).

For (2), assume that $B$ is rational and $K_{\tilde{Y}}+\Theta$
has a Zariski decomposition relative to $S$. In particular, 
$\tilde{B}$ has rational coefficients. By Theorem~\ref{mr} 
applied to $(\tilde{X}/S,\tilde{B})$, we infer that $\Theta$ is 
rational and the $\cO_S$-algebra 
$
\bigoplus_{n=0}^\infty (\pi\circ\mu)_*\cO_{\tilde{Y}}
(nK_{\tilde{Y}}+n\Theta)
$
is finitely generated.
\end{proof}

\begin{thm}\label{pef}
Let $X$ be a nonsingular variety, endowed with an $\R$-divisor 
$B$ such that $\Supp(B)$ is a simple normal crossings divisor. 
Assume that $Y$ is a component of $B$ of multiplicity one 
and $\lfloor B-Y\rfloor=0$. Denote $B_Y=(B-Y)\vert_Y$, so 
that by adjunction we have
$
(K+B)\vert_Y=K_Y+B_Y.
$
Let $\pi\colon X\to S$ be a projective surjective morphism,
$H$ a $\pi$-free divisor on $X$ and $n\ge 2$ an integer such 
that $\{n B\}\le B$. Assume that
\begin{itemize}
\item[(a)] $nK+nB+H\sim_\Q A+C$, where $A$ is 
a $\pi$-ample $\R$-divisor and $C$ is an 
effective $\R$-divisor with $\mult_Y(A)=\mult_Y(C)=0$.
\item[(b)] $(Y,B_Y)$ has canonical singularities 
in codimension at least two.
\item[(c)] $n\ge (1-b)^{-1}$, where $b$ is the maximum
multiplicity of the components of $B_Y$.
\end{itemize}
Define $\Theta_n=\max(B_Y-\lim_{i\to \infty}
\frac{(\bFix(i(K+B+\frac{H}{n}))\vert_Y)_Y}{i},0)$.
Then 
$$
\im(\pi_*\cO_X(nK+nB+H)\to \pi_*\cO_Y(nK_Y+nB_Y+H\vert_Y))
=
\pi_*\cO_Y(nK_Y+n\Theta_n+H\vert_Y).
$$
\end{thm}

\begin{proof} We may assume that $S$ is affine. Let 
$D\in \vert H\vert$ be a general member such that 
$D\vert_Y$ does not contain any prime component of
the $\R$-divisor
$
\lim_{i\to \infty}\frac{(\bFix(iK+iB)\vert_Y)_Y}{i}.
$
Define $B_n=B+\frac{1}{n}D$ and $B_{n,Y}=(B_n-Y)\vert_Y$.
By (c) and since $H$ is $\pi$-free, $(Y,B_{n,Y})$ is also 
canonical in codimension a least two. By construction, 
$B_{n,Y}=B_Y+\frac{1}{n}D\vert_Y$. By the choice of $D$,
we have 
$$
\max(B_{n,Y}-
\lim_{i\to \infty}\frac{(\bFix(i(K+B+\frac{D}{n}))\vert_Y)_Y}{i},0)
=\Theta_n+\frac{1}{n}D\vert_Y.
$$
Now apply Theorem~\ref{mr} to $(X,B_n)$ and $n$.
\end{proof}

%%%%%%%%%%%%%%%%%%%%%%%%%%%%%%%%%%%%%%%
%%%%%%%%%%%%%%%%%%%%%%%%%%%%%%%%%%%%%%%

\end{document}